\documentstyle[12pt]{article}
\newcommand{\setleftmargin}[1]{
	\addtolength{\textwidth}{\oddsidemargin}
	\addtolength{\textwidth}{1in}
	\addtolength{\textwidth}{-#1}
	\setlength{\oddsidemargin}{-1in}
	\addtolength{\oddsidemargin}{#1}
	\setlength{\evensidemargin}{\oddsidemargin}
}
\newcommand{\setrightmargin}[1]{
	\setlength{\textwidth}{8.5in}
	\addtolength{\textwidth}{-\oddsidemargin}
	\addtolength{\textwidth}{-1in}
	\addtolength{\textwidth}{-#1}
}
\newcommand{\settopmargin}[1]{
	\addtolength{\textheight}{\topmargin}
	\addtolength{\textheight}{1in}
	\addtolength{\textheight}{\headheight}
	\addtolength{\textheight}{\headsep}
	\addtolength{\textheight}{-#1}
	\setlength{\topmargin}{-1in}
	\addtolength{\topmargin}{-\headheight}
	\addtolength{\topmargin}{-\headsep}
	\addtolength{\topmargin}{#1}
}
\newcommand{\setbottommargin}[1]{
	\setlength{\textheight}{11in}
	\addtolength{\textheight}{-\topmargin}
	\addtolength{\textheight}{-1in}
	\addtolength{\textheight}{-\footskip}
	\addtolength{\textheight}{-#1}
}
\newcommand{\setallmargins}[1]{
	\settopmargin{#1}
	\setbottommargin{#1}
	\setleftmargin{#1}
	\setrightmargin{#1}
}
\setallmargins{1.2in}
\title{ Simple Lie algebras which generalize Witt algebras
}
\author{ Ki-Bong Nam 
     \thanks{ 
    Department of Mathematics,
    University of Wisconsin, Madison, WI 53706}}
\begin{document}
\maketitle
\begin{abstract}
We introduce a new class of simple
Lie algebras $W(n,m)$ (see Definition 1) 
that generalize the Witt algebra by using "exponential"
functions,
and also a subalgebra $W^*(n,m)$ thereof;
and we show each derivation of $W^*(1,0)$ can be written as a sum
of an inner derivation and a scalar derivation (Theorem. 2) \cite{Nam}.
The Lie algebra $W(n,m)$ is $Z$-graded and is infinite growth \cite{Kac}.
\end{abstract}

\newtheorem{lemma}{Lemma}
\newtheorem{prop}{Proposition}
\newtheorem{thm}{Theorem}
\newtheorem{coro}{Corollary}
\newtheorem{definition}{Definition}

\section{Introduction}

Let $F$ be a field of characteristic zero.
Let us start by 
recalling the ${\bf general }$ ${\bf algebra}$ $W^+(n)$. 
(For more details, please
refer to \cite{Rud}.)

The Lie algebra $W(n)$ has a basis
$$\{x_1^{i_1}\cdots x_n^{i_n} \partial_{k}\quad |\quad i_1,\cdots,i_n\in Z,
1\leq k\leq n\},$$
and a Lie bracket on basis elements given by
\begin{eqnarray}\label{s1}
& &[x_1^{i_1}\cdots x_n^{i_n} \partial_{k},
x_1^{j_1}\cdots x_n^{j_n} \partial_{t}] \\ \nonumber
&=&j_k x_1^{i_1+j_1}\cdots  x_k^{i_k+j_k-1} x_{k+1}^{i_{k+1}+j_{k+1}}\cdots x_n^{i_n+j_n} \partial_{t}\\ \nonumber
&-&i_t x_1^{i_1+j_1}\cdots  x_t^{i_t+j_t-1} x_{t+1}^{i_{t+1}+j_{t+1}}\cdots x_n^{i_n+j_n} \partial_{k} 
\end{eqnarray}

%for any pair of basis elements
%$$x_1^{i_1}\cdots x_n^{i_n} \partial_{k},
%x_1^{j_1}\cdots x_n^{j_n} \partial_{t} \in W(n),$$
\noindent
where $\partial_t$ is the partial derivative with
respect to $x_t,$ $1\leq t \leq n.$

%%%%%%%%%%%%%%%%%%%%%%%%%%%%%%%%%%%%%%%%%%%%%%%%%%%%%%%%%%%%%%%%%%%%%%%
It is well known that the Witt algebra $W^+(n)$ 
is 
defined to be 
%(which is just Witt algebra)
the subalgebra of $W(n)$,
with the basis
$$\{x_1^{i_1}\cdots x_n^{i_n} \partial_{k}\quad |\quad i_1,\cdots,i_n\in Z^+\cup \{0\},
1\leq k \leq n\}.$$
%and Lie brackets on basis elements given by

%\begin{eqnarray}\label{s4}
%& &[x_1^{i_1}\cdots x_n^{i_n} \partial_{k},
%x_1^{j_1}\cdots x_n^{j_n} \partial_{t}] \\ \nonumber
%&=&j_k x_1^{i_1+j_1}\cdots  x_k^{i_k+j_k-1} x_{k+1}^{i_{k+1}+j_{k+1}}\cdots x_n^{i_n+j_n} \partial_{t}\\ \nonumber
%&-&i_t x_1^{i_1+j_1}\cdots  x_t^{i_t+j_t-1} x_{t+1}^{i_{t+1}+j_{t+1}}\cdots x_n^{i_n+j_n} \partial_{t} 
%\end{eqnarray}
%where $\partial_t$ is the partial derivation with
%respect to $x_t.$ 

In \cite{Kaw} Kawamoto generalizes Witt algebras; we
call the algebra $W(G,I)$ in \cite{Kaw} 
Kawamoto algebra.
Dokovic and Zhao \cite{Do} show that their algebras generalize Kawamoto
algebras. But their algebras are different from our algebras
because their algebras have ad semi-simple elements
(see \cite{Do}), whereas our algebras have none
.

The main results of this paper are following:

\noindent
{\bf Theorem 1} The Lie algebras
$W(n,m), W(n,m,r,s)$ and $W^*(n,m)$ are simple 
Lie algebras.

\bigskip

\noindent
(See Definition 1, Definition 2 and Definition 3.)

\bigskip

\noindent
{\bf Theorem 3} 
Each derivation of $W^*(1,0)$ is a sum of
an inner derivation and a scalar derivation.

\section{Generalized Witt algebras}

We now give the definition of the algebra $W(n,m).$
\begin{definition}
For any pair of natural numbers $n,m$, 
define the Lie algebra 
$W(n,m)$ to be the algebra spanned linearly by a basis 
\begin{eqnarray*}
& &\{e^{a_1 x_1} \cdots e^{a_n x_n}x_1^{b_1}\cdots x_{m+n}^{b_{m+n}} 
\partial_{i}
\quad |\quad a_1,\cdots, b_1,\cdots,b_{m+n}\in Z, 1\leq i \leq m+n \},
\end{eqnarray*}
with Lie bracket on basis elements given by
\begin{eqnarray}\label{s5}
& &[e^{a_1 x_1} \cdots e^{a_n x_n}x_1^{l_1}\cdots x_{m+n}^{l_{m+n}} 
\partial_{i} \quad ,
\quad e^{b_1 x_1} \cdots e^{b_n x_n}x_1^{t_1}\cdots x_{m+n}^{t_{m+n}} 
\partial_{j}] \\ \nonumber
&=&b_i e^{a_1x_1+b_1 x_1} \cdots e^{a_n x_n +b_n x_n}x_1^{l_1+t_1}\cdots x_{m+n}^{t_{m+n}+l_{n+m} }
\partial_{j} \\ \nonumber
&+&t_i e^{a_1x_1+b_1 x_1} \cdots e^{a_n x_n +b_n x_n}x_1^{l_1+t_1}\cdots x_{m+n}^{t_{m+n}+l_{n+m} }
x_i^{-1}\partial_{j} \\ \nonumber
&-&a_j e^{a_1x_1+b_1 x_1} \cdots e^{a_n x_n +b_n x_n}x_1^{l_1+t_1}\cdots x_{m+n}^{t_{m+n}+l_{n+m} }
\partial_{i} \\ \nonumber
&-&l_i e^{a_1x_1+b_1 x_1} \cdots e^{a_n x_n +b_n x_n}x_1^{l_1+t_1}\cdots x_{m+n}^{t_{m+n}+l_{n+m} }
x_j^{-1}\partial_{i}, 
\end{eqnarray}
where $b_i=0$ if $n+1\leq i \leq n+m,$
and $a_j=0$ if $n+1\leq j \leq n+m$ (1).
%%%%%%%%%%%%%%%%%%%%%%%%%%%%%%%%%%%%%%%
\end{definition}

\begin{definition}
Define $W^*(n,m)$ 
to be the Lie subalgebra of the Lie
algebra $W(n,m)$ 
%as the Lie sub-algebra 
spanned by the basis
\begin{eqnarray}\label{s7}
& & e^{a_1 x_1} 
\cdots e^{a_n x_n}x_1^{b_1}\cdots x_{m+n}^{b_{m+n}} \partial_{i} 
\end{eqnarray}
where $a_1,\cdots, a_n\in Z, b_{1},\cdots,b_{m+n} \in N, \hbox { and } 1\leq i \leq m+n$.
\end{definition}

\noindent
Remarks: It can be checked that the Lie bracket of
any two such basis elements produces a linear combination
of the same form.
\begin{definition}
Define 
$W(n,m,r,s)$  
as the Lie subalgebra 
of $W(n,m)$  
spanned by the basis
\begin{eqnarray*}
& &\{e^{a_1 x_1} \cdots e^{a_n x_n}x_1^{i_1}\cdots x_r^{i_r} x_{r+1}^{i_{r+1}} \cdots x_n^{i_n} x_{n+1}^{i_{n+1}} \cdots x_{n+s}^{i_{n+s}}x_{n+s+1}^{i_{n+s+1}}\cdots x_{m+n}^{b_{m+n}} 
\partial_{k} \quad |\\
& &\quad a_1,\cdots, a_{n}\in Z, 
i_1,\cdots, i_{r},i_{n+1},\cdots ,i_{n+s}\in Z, 
i_{r+1},\cdots, i_{n},i_{n+1+s},\cdots ,i_{n+s}\in N, \\
& &1\leq k\leq  m+n\}
\end{eqnarray*}
where $N$ is the set of non-negative integers and 
$r\leq n\in N, s\leq m\in N.$
\end{definition}

%\begin{definition}
%Define 
$W(0,s+r,0,r)$  
is the Lie subalgebra 
of $W(n,m)$  
spanned by the basis
\begin{eqnarray*}
& &\{x_1^{i_1}\cdots x_r^{i_r} x_{r+1}^{i_{r+1}} \cdots x_{r+s}^{i_{r+s}} 
\partial_{j}\quad | \\
& &\quad i_1,\cdots, i_{r}\in Z, 
i_{r+1},\cdots, i_{r+s}\in N, 
1\leq j \leq r+s \},
\end{eqnarray*}
where $r, s \in N.$
%\end{definition}

Note that the Lie algebra $W(0,s,0,r)$ is 
isomorphic to a Lie algebra in the class $W^*$ of 
\cite{Os}.

%%%%%%%%%%%%%%%%%%%%%%%%%%%%%%%%%%%%%%%%%%%%%%%%%%%
%%%%%%%%%%%%%%%%%%%%%%%%%%%%%%%%%%%%%%%%%%%%%%%%%%%
%It is easy to verify that the linear subspace generated by $(1)$
%is a Lie subalgebra of $W(n,m).$

\noindent
%Let us denote $Z^k$ as $k-$times 
%$Z \hbox { x } Z \hbox { x } \cdots \hbox { x } Z$ of $Z.$
Let $Z^k$ denotes the product of $k$ factors of $Z.$
We consider an $n-$tuple
$\alpha =(a_1,\cdots ,a_n)\in Z^n$ 
%\hbox { x } Z \hbox { x } \cdots \hbox { x } Z$
and an 
($n+m$)$-$tuple
%$\hbox { ( }n+m \hbox { ) }-$tuple
$\beta =(b_1,\cdots ,b_{n+m})\in Z^{n+m}.$ 
%\hbox { x } Z \hbox { x } \cdots \hbox { x } Z.$
If we put 
$$e^{\alpha}:=
e^{a_1 x_1} \cdots e^{a_n x_n}$$
and
 $$x^{\beta}:=
x_1^{b_1}\cdots x_{m+n}^{b_{m+n}},$$
then any element $l$ in $W(n,m)$
has the form $$l=\sum _{\alpha,\beta,p}a_{\alpha,\beta,p}e^{\alpha} x^{\beta}\partial_p,$$
for $a_{\alpha,\beta,p}\in F$.
Let us write $e^{ax_i}\partial_p=e^{a\epsilon_i}\partial_p,$
and
$x_i^b\partial_p=x^{b\epsilon_i}\partial_p,$
for $1\leq p\leq n+m,$ where $a,b\in Z.$

The Lie algebra $W(n,m)$ has a $Z^k-$gradation 
for $1\leq k \leq n $ as follows: 
\begin{eqnarray}\label{s9}
& &W(n,m)=\bigoplus_{(a_1,\cdots ,a_k
)\in Z^{k}} 
W_{(a_1,\cdots ,a_k
%,a_{k+1},\cdots ,a_{n}
)}  
\end{eqnarray}
where
$W_{(a_1,\cdots ,a_k
%,a_{k+1},\cdots ,a_{n}
)} $ is a subspace of $W(n,m)$ with basis
\begin{eqnarray*}
%& &W_{(a_1,\cdots ,a_k
%,a_{k+1},\cdots ,a_{n}
%)\in Z^{k}} 
& &\{ e^{\alpha x}x^{\beta}\partial_k 
\quad | \quad
\alpha =(a_1,\cdots ,a_k,a_{k+1},\cdots ,a_n)\in Z^n, \beta \in
Z^{n+m}, \\
& &1\leq k\leq n+m \}.
\end{eqnarray*}

\noindent

\noindent
%Note that $W_0$ is isomorphic to the Witt algebra $W(n+m).$
Let us call $W_{\alpha}$ the $\alpha$-homogeneous component
of $W(n,m)$ and
elements in $W_{\alpha}$ the
$\alpha$-homogeneous elements. 
Note that the $(0,\cdots ,0)$-homogeneous component
is isomorphic to the Witt algebra $W(n+m).$
From now on let us denote the $(0,\cdots ,0)$-homogeneous component
as the $0$-homogeneous component.

%%%%%%%%%%%%%%\begin{definition}\label{s15}
%%%%%%%%%%%%%For any two basis elements
%%%%%%%%%%%%%$l_1\in W_{(a_1,\cdots ,a_n)}$
%%%%%%%%%%%%%and $l_2\in W_{(b_1,\cdots ,b_n)}$ we define
%%%%%%%%%%%%%an lexicographic ordering $>_o$ as follows:
%%%%%%%%%%%%%$$\hbox { if } l_1=e^{a_1 x_1}\cdots e^{a_n x_n}x_1^{l_1}\cdots x_{m+n}^{l_{n+m}}\partial_k \hbox { and }
%%%%%%%%%%%%%l_2=e^{b_1 x_1}\cdots e^{b_n x_n}x_1^{t_1}\cdots x_{m+n}^{t_{n+m}}\partial_p,$$
%%%%%%%%%%%%%\begin{eqnarray}
%%%%%%%%%%%%%& &\hbox { then } l_1>_o l_2 \hbox { if and only if }\\ \nonumber
%%%%%%%%%%%%%& &a_1>b_1, \hbox { or } a_1=b_1 \hbox { and } a_2>b_2, \hbox { or } \ldots ,\hbox { or } \\ \nonumber
%%%%%%%%%%%%%& &l_{n+m-2}=t_{n+m-2} \hbox { and }  l_{n+m-1}>t_{n+m-1}, \hbox { or } \\ \nonumber
%%%%%%%%%%%%%& &l_{n+m-1}=t_{n+m-1} \hbox { and }  l_{n+m}>t_{n+m}, \hbox { or } l_{n+m}=t_{n+m} \hbox { and } k<p. 
%%%%%%%%%%%%%\end{eqnarray}
%%%%%%%%%%%%%\end{definition}
We introduce a lexicographic ordering of the basis
elements of $W(n,m)$ as follows:
\begin{eqnarray}\label{s15}
& &\hbox { For } l_1=e^{a_1 x_1}\cdots e^{a_n x_n}
x_1^{b_1}\cdots x_{m+n}^{b_{n+m}}\partial_p \\ \nonumber
& &\hbox { and }
 l_2=e^{c_1 x_1}\cdots e^{c_n x_n}
x_1^{d_1}\cdots x_{m+n}^{d_{n+m}}\partial_q,  \\ \nonumber
& &\hbox { we have } l_1\geq l_2 \\ \nonumber
& &\hbox { if } (a_1,\cdots ,a_n,b_1,\cdots ,b_{n+m},p) \geq
(c_1,\cdots ,c_n,d_1,\cdots ,d_{n+m},q)
\end{eqnarray}
by the natural lexicographic ordering in
$Z^n\times Z^{n+m}\times Z.$

\noindent
%%%%%%%%%%%%%Using the above ordering on gradations, for any
%%%%%%%%%%%%%$e^{\alpha x} x^{\beta}\partial_i,
%%%%%%%%%%%%%e^{\phi x} x^{\psi}\partial_j \in 
%%%%%%%%%%%%%W(n,m),$
%%%%%%%%%%%%%we can induce lexicographic ordering
%%%%%%%%%%%%%on $W(n,m) $
%%%%%%%%%%%%%%\alpha >_o \phi$ or $\beta >_o\psi$ etc., using
%%%%%%%%%%%%%using the lexicographic ordering 
%%%%%%%%%%%%%on the coordinates of $\alpha =(a_1,\cdots ,a_n),$
%%%%%%%%%%%%%$\phi =(b_1,\cdots ,b_n)\in Z^n$
%%%%%%%%%%%%%etc.
%%%%%%%%%%%%%
\noindent
For any element $l\in W(n,m)$, $l$ can be written as follows using
the ordering and gradation,
\begin{eqnarray}\label{s20}
& &l=\sum_{\beta _{1 }} a_{\alpha_{1},\beta_{1},p} e^{\alpha_{1}}x^{\beta_{1}}\partial_p\\ \nonumber
&+&\sum_{\beta_{2}} a_{\alpha_{2},\beta_{2},p} e^{\alpha_{2}}x^{\beta_{2}}\partial_r \\ \nonumber
&+& \cdots \\ \nonumber
&+&\sum_{\beta_{s}} a_{\alpha_{s},\beta_{s},r} e^{\alpha_{s}}x^{\beta_{s}}\partial_q, 
\end{eqnarray}
where
$\alpha_{1}>\alpha_{2}>\cdots >\alpha_{s}.$ 

\noindent
%Let us define the {\bf string number} $st(l)$ of $l$
%as the cardinality of
%$\{\alpha_{1i},\alpha_{2s},\cdots ,\alpha_{sq}\}.$ 
Let us define the {\bf string number} $st(l)=s$ for $l.$

%\bigskip
%\noindent
%{\bf Example}
%For $l=e^{3x_1} x_5^5\partial_3 +e^{3x_1}x_3^{2}\partial_7
%+e^{4x_3} \partial_2\in W(n,m),$ we have $st(l)=2.$

%\bigskip
\noindent
For any $l\in W(n,m)$, let us define a largest
power $lp_{(a_1,\cdots ,a_n)}(l)$ as the largest power of polynomials
in $(a_1,\cdots ,a_n)$-homogeneous elements.

\bigskip
\noindent
Similarily for any $l\in W(n,m)$, we can define a largest
power $lp_l(l)$ as the largest power of polynomials
in $l$-homogeneous elements.
%%%%%`For any $l\in W(n,m)$, let us define the largest
%%%%%`power $lp_{\alpha_i}(l)$ as the largest power
%%%%%`of $x_1,\cdots ,x_{n+m}$ in $\alpha_i$-homogeneous
%%%%%`component of $l$, and $lp(l)=\max_{\alpha_i}

%\bigskip
%noindent
%{\bf Example}
%For $l=8 e^{6x_1}x_1^5 x_2^7\partial_1 
%+9 e^{7x_2}x_2^7 x_2^{-8}\partial_4 
%+6 x_1^3 x_5^{-2}\partial_3\in W(n,m)$, 
%we have $lp_{(0,\cdots ,0)}(l)=3$ and $lp_l(l)=7.$

\bigskip

%%%%%%%%%%%%%%%%%%%%%\begin{lemma}
%%%%%%%%%%%%%%%%%%%If $l=\sum_{\alpha,\beta,p} a_{\alpha,\beta,p}e^{\alpha}x^{\beta}\partial_p$
%%%%%%%%%%%%%%%%%%%is an element in $W(n,m)$ for 
%%%%%%%%%%%%%%%%%%%$e^{\alpha}x^{\beta}\partial_p\in W_{\alpha}$ and for any
%%%%%%%%%%%%%%%%%%%$e^{-a_s \epsilon_s}x^{\beta'}\partial_t\in W(n,m)$ 
%%%%%%%%%%%%%%%%%%%with $1\leq t\leq n+m$ and $i+1 \leq s \leq n,$
%%%%%%%%%%%%%%%%%%%then
%%%%%%%%%%%%%%%%%%%$$[e^{-a_s \epsilon_s}x^{\beta'}\partial_t,l] 
%%%%%%%%%%%%%%%%%%%=\sum_{\alpha',\beta',s} a_{\alpha',\beta',p}e^{\alpha'}x^{\beta'}\partial_p$$
%%%%%%%%%%%%%%%%%%%where $ a_{\alpha',\beta',s} \in F$
%%%%%%%%%%%%%%%%%%%and 
%%%%%%%%%%%%%%%%%%%$e^{\alpha'}x^{\beta'}\partial_p\in W_{\alpha}.$
%%%%%%%%%%%%%%%%%%%\end{lemma}
%%%%%%%%%%%%%%%%%%%{\it Proof.}
%%%%%%%%%%%%%%%%%%%This is clear by observing the exponential parts of
%%%%%%%%%%%%%%%%%%%$[e^{-a_s \epsilon_s}x^{\beta'}\partial_t,l]$ 
%%%%%%%%%%%%%%%%%%%.
%%%%%%%%%%%%%%%%%%%\quad $\Box$
%%%%%%%%%%%%%%%%%%%
\noindent
%Let us omit the subscripts from here, 
For the sake of clarity,
from now on we generally omit subsripts.
%%%%%%%%%%%%%%%%%%%%%%%%%
\begin{lemma}
If $l\in W(n,m)$ is any non-zero element
then the ideal $<l>$ 
generated by $l$ contains an element $l'$ 
whose powers of polynomial parts
are positive integers.
\end{lemma}
{\it Proof.}
%%%%%%%%%%%%%%%%%Since $l$ is a non-zero element of $W(n,m)$, without loss of generality
%%%%%%%%%%%%%%%%%we may assume that there is a non-zero 
%%%%%%%%%%%%%%%%%%%%%$(0,\cdots ,a_j,\cdots ,a_n)$-homogeneous element.
%%%%%%%%%%%%%%%%We can take an element $s=x^{u} \partial_i$
%%%%%%%%%%%%%%%%%with large positive integers $u\in W^{n+m}$. 
%%%%%%%%%%%%%%%%%Then $$[l,s]\neq 0,$$ is the required the element
%%%%%%%%%%%%%%%%%in the lemma.
Since $l$ is a non-zero element of $W(n,m)$, there
is a non-zero $\alpha$-homogeneous component.
Then we can take an element $s=x^{u}\partial_i$
with large positive integers $u\in Z^{n+m}$
such that 
$l'=[s,l]\neq 0.$
\quad $\Box$

%Let us call the following lemma as string lemma.
\begin{lemma}
%%%%%%%%%%%%%%%%$i\in \{1,\cdots, n+m\}$ is $W(n,m).$
If an ideal $I$ of $W(n,m)$ contains
any $\partial_i$ where $(1\leq i\leq n+m)$, then 
$I=W(n,m)$.
\end{lemma}
{\it Proof.}
Since the Witt algebra $W(n+m)$ is simple \cite{Kaw},\cite{Os},\cite{Rud}
and $\partial_i\in I,$
the Witt algebra $W(n+m)\subset I.$
For any basis element $e^{\alpha}x^{\beta}\partial_t$ of $W(n,m),$
we can assume $a_k \neq 0$ for some $k\in \{1,\cdots,n\}$
because if all $a_k=0,$ then the basis $x^{\beta}  \partial_t \in W(n+m)$ 
where $W(n+m)$ is a simple Lie algebra \cite{Kaw}, \cite{Os}.

\noindent
If $\beta=(0,\cdots ,0)$, then there is $a_k\neq 0.$
We have
$$[\partial_k\quad ,\quad e^{\alpha} \partial_t]=a_k e^{\alpha}\partial_t\in I.$$

For any basis element $e^{\alpha}x^{\beta}\partial_t$
such that $a_u\neq 0$ for $1\leq u \leq n,$ we have
\begin{eqnarray}\label{s23}
& &[\partial_u, e^{\alpha}x^{\beta}\partial_t]=
a_u e^{\alpha}x^{\beta}\partial_t
+b_u e^{\alpha}x^{\beta-\epsilon_u}\partial_t.
\end{eqnarray}
\begin{eqnarray}\label{s24}
& &[e^{a_u\epsilon_u}\partial_u, e^{\alpha-a_u\epsilon_u}x^{\beta}\partial_t]=
b_u e^{\alpha}x^{\beta-\epsilon_u}\partial_t
-\delta_{u,t}a_t e^{\alpha}x^{\beta}\partial_u
\end{eqnarray}
where $\delta_{u,t}$ is the Knonecker delta.

\noindent
Case I.
If $u\neq t,$
$$(\ref{s23})-(\ref{s24})=
a_ue^{\alpha}x^{\beta}\partial_t\in I.$$

\noindent
Case II.
If $u= t,$
$$(\ref{s23})-(\ref{s24})=
2a_ue^{\alpha}x^{\beta}\partial_t\in I.$$

%\noindent
%Case II. If $\beta_k\neq 0$ for some $k\in \{1,\cdots,n\}$ and $\beta_k=1$, then
%$$[\partial_k\quad ,\quad e^{\alpha} x^{\beta}\partial_t]=\alpha_k e^{\alpha}x^{\beta}\partial_t +e^{\alpha}x^{\beta -\epsilon_k}\partial_t\in I.$$
%Thus $a_k e^{\alpha}x^{\beta}\partial_t\in I.$
%
%
%\noindent
%Similarly we have arguing inductively
%$e^{\alpha}x^{\beta +s\epsilon_k}\partial_t\in I,
%$ for $s\in N.$
%Using $x^{\epsilon_k}\partial_k$ instead of $\partial_k$, 
%we have $e^{\alpha}x^{\beta}\partial_t\in I$ where
%$a_k\neq 0$ and $\beta_k\in Z-N.$ 
Therefore, we have proved the lemma.
\quad $\Box$

%Since that the Witt algebra $W(n+m)$ is simple \cite{Kaw}, \cite{Os}
%and $\partial_i\in I$
%,the Witt algebra $W(n+m)\subset I.$
%For any basis element $e^{\alpha}x^{\beta}\partial_t$ of $W(n,m)$,
%we can assume $\alpha_k \neq 0$ for some $k\in \{1,\cdots,n\}$
%because if all $a_k=0,$ then the basis $x^{\beta}  \partial_t \in W(n+m)$ 
%where $W(n+m)$ is a simple Lie algebra \cite{Kaw}, \cite{Os}.
%
%
%\noindent
%Case I. If $\beta_k=0$ for all $k$, then
%$$[\partial_k\quad ,\quad e^{\alpha} x^{\beta}\partial_t]=\alpha_i e^{\alpha}x^{\beta}\partial_t\in I.$$
%
%
%\noindent
%Case II. If $\beta_k\neq 0$ for  some $k\in \{1,\cdots,n\}$ and take $\beta_k=1$, then
%$$[\partial_k\quad ,\quad e^{\alpha} x^{\beta}\partial_t]=\alpha_k e^{\alpha}x^{\beta}\partial_t +e^{\alpha}x^{\beta -\epsilon_k}\partial_t\in I.$$
%Thus $e^{\alpha}x^{\beta}\partial_t\in I.$
%
%
%\noindent
%Similarly we have $e^{\alpha}x^{\beta +s\epsilon_k}\partial_t\in I,
%$ for $s\in N.$
%Using $x^{\epsilon_k}\partial_k$ instead of $\partial_k$, 
%we have $e^{\alpha}x^{\beta}\partial_t\in I$ where
%$\alpha_k\neq 0$ and $\beta_k\in Z-N.$ 
%Therefore, we have proved the lemma.
\quad $\Box$

\begin{lemma}
%%%%%%%%%%%%%%%%%%%%%%%%%%%%%%For a non-zero element $\sum c_{\alpha, \beta}e^{\alpha x}x^{\beta }\partial_j \in
%%%%%%%%%%%%%%%%%%%%%%%%%%%%%%%%%%%%%%%%%%%%%%%%%%%%%%%%%%%%W(n,m)$ where \\ 
%%%%%%%%%%%%%%%%%%%%%%%%%%%%%%$\alpha =(0,\cdots ,0,a_u,a_{u+1},\cdots, a_n)$  
%%%%%%%%%%%%%%%%%%%%%%%%%%%%%%with $a_u\neq 0$,then 
%%%%%%%%%%%%%%%%%%%%%%%%%%%%%%$[\partial_u,\sum c_{\alpha, \beta}e^{\alpha x}x^{\beta }\partial_j ]\neq 0.$
If $\sum c_{\alpha, \beta}e^{\alpha x}x^{\beta }\partial_j $
is a non-zero element of $W(n,m)$ where \\
$\alpha =(a_1,\cdots ,a_{u},\cdots, a_n)$  with
$a_u\neq 0$, then
$[\partial_u \quad ,\quad \sum c_{\alpha, \beta}e^{\alpha x}x^{\beta }\partial_j ]\neq 0.$
\end{lemma}
{\it Proof.}
We have
\begin{eqnarray}\label{s25}
& &[\partial_u\quad ,\quad \sum c_{\alpha, \beta}e^{\alpha }x^{\beta }\partial_j] \\ \nonumber
&=&\sum a_u c_{\alpha, \beta}e^{\alpha }x^{\beta }\partial_j 
+\sum b_u c_{\alpha, \beta}e^{\alpha }x^{\beta -\epsilon_u}\partial_j\neq 0, 
\end{eqnarray}
%%%%%%%%%%%%%%%%%since $c_{\alpha,\beta}a_u\neq 0$ and the maximality of the
%%%%%%%%%%%%%%%%%%%$\alpha$-homogeneous element in (7).
by $c_{\alpha,\beta}a_u\neq 0$ and the 
ordering in $W(n,m)$.
\quad $\Box$

\begin{thm}
The Lie algebras $W(n,m),$ $W(n,m,r,s)$ and $W^*(n,m)$
are simple Lie algebras.
\end{thm}
{\it Proof.}
%Let us prove this theorem in three steps.
If $I$ is a non-zero ideal of $W(n,m),$  
(or $W^*(n,m)$ or $W(n,m,r,s)$) then by Lemma 1 there is
a non-zero element $l$ in $I$ such that all the polynomial
parts are positive powers.
%We prove this theorem by induction
%on $st(l).$
We show that $I=W(n,m)$ by induction on
$st(l).$

\noindent
%{\bf Step 1}
%%%%%%%%%%%%%%%%%%%%%%%%%%%%%%Assume $st(l)=1$ for any $l\in I$, then the theorem is true.

\noindent
%{\it Proof.}
If $st(l)=1,$ then $l$ can be written as follows:
$$l=c_1 e^{\alpha }x^{\beta_1}\partial_{j_1}+\cdots + c_te^{\alpha }x^{\beta_t}\partial_{j_t},$$
where $c_1,\cdots ,c_t\in F$
$\alpha =(a_1,\cdots ,a_n) \in Z^n,$
and $\beta_1,\cdots ,\beta_{m+n} \in Z^{n+m}.$
If $\alpha =(0,\cdots ,0),$ then
$l\in W_0^n = W(n+m)$
and $I \cap W(n+m)$ is a non-zero ideal of $W(n+m).$
Thus, $W(n+m)\subset I$ 
%the theorem is proved by lemma 3, and the
by the simplicity of $W(n+m)$
and $I=W(n,m)$ by Lemma 2.

\noindent
%%%%%%%%%%%%%%%%%%%%$If $\alpha \neq (0,\cdots ,0),$ then
If $\alpha =(a_1,\cdots ,a_n) \neq (0,\cdots ,0),$ then
without loss of generality  we can assume
%%%%%%%%%%%%%%%%%%%%$\alpha =(a_1,a_2,\cdots ,a_n)$ such that 
$a_1\neq 0.$
Then,
\begin{eqnarray}\label{s30}
& &l_1=[e^{-\alpha }\partial_1 \quad ,
\quad c_1 e^{\alpha }x^{\beta_1}\partial_{j_1}+\cdots +c_te^{\alpha }x^{\beta_t}\partial_{j_t}]\\ \nonumber
&=&c_1 a_{j_1}x^{\beta_1}\partial_{1}
+c_1 a_{1}x^{\beta_1}\partial_{j_1}
+c_1 b_{11}x^{\beta_1-\epsilon_1}\partial_{j_1}\\ \nonumber
&+&\cdots \\ \nonumber
&+&c_1 a_{j_t}x^{\beta_t}\partial_{1}
+c_1 a_{1}x^{\beta_t}\partial_{j_t}.
+c_1 b_{t1}x^{\beta_t-\epsilon_1}\partial_{j_t}.
\end{eqnarray}
%%%%%%%%%%%%%%%%%%%%%%%If $j_1=1$ and $a_{j_1}\neq 0$, 
If $a_{j_1}\neq 0$ and $j_1= 1$, 
then we have from (\ref{s30}),
\begin{eqnarray} \label{s35}
& &l_1=2a_1 c_1 x^{\beta_1}\partial_1
+c_2 a_{j_2}x^{\beta_2}\partial_1 \\ \nonumber
&+&\cdots \\ \nonumber
&+&c_t a_{j_t}x^{\beta_t}\partial_1
+c_1 b_{11}x^{\beta_1-\epsilon_1}\partial_{1}\\ \nonumber
&+&c_2 b_{21}x^{\beta_t-\epsilon_1}\partial_{j_2} +\cdots +c_t b_{t1} x^{\beta _t-\epsilon_1}\partial_{j_t} \neq 0,
\end{eqnarray}
since $2a_1c_1\neq 0$ and the maximality of $x^{\beta_1}\partial_1.$
%%%%%%%%%%%%%%%%%%%%%Since $l_1\in W(n+m)$, the theorem
%%%%%%%%%%%%%%%%%%%%%is proved again by lemma 3 and the simplicity of $W(n+m)$ 
%%%%%%%%%%%%%%%%%%%%%.
Since $l_1\in W(n+m)$, we have $I=W(n,m)$
by the simplicity of $W(n+m)$ 
and by Lemma 2.  

\noindent
If $a_{j_1} \neq 0$
and $j_1\neq 1,$ then from (\ref{s30}),
\begin{eqnarray}\label{s40}
& &l_1=c_1 a_{j_1} x^{\beta_1}\partial_1
+c_1 a_{1}x^{\beta_1}\partial_{j_1} 
+c_1 b_{1}x^{\beta_1-\epsilon_1}\partial_{j_1} 
\\ \nonumber
&+&\cdots \\ \nonumber
&+&c_t a_{j_t}x^{\beta_t}\partial_1
+c_t a_{1} x^{\beta_t}\partial_{j_t}
+c_t b_{t1} x^{\beta_t-\epsilon_1}\partial_{j_t}.
%&+&\cdots \\ \nonumber
%&+&c_t a_{1} x^{\beta_t}\partial_{j_t}\\ .
\end{eqnarray}
%%%%%%%%%%%%%%%%%%%Since $l_1\neq 0$, $x^{\beta_1}>x^{\beta_2}$ by assumption.
%%%%%%%%%%%%%%%%%%Thus, we have proved the theorem by lemma 3 and the simplicity
%%%%%%%%%%%%%%%%%%%%%of $W(n+m).$
If
 $a_{j_1} \neq 0$ and $j_1\neq 1,$ then
%from (\ref{30}), 
$c_1 a_{j_1}x^{\beta_1}\partial_1$
is the non-zero maximal element by
the ordering in $l_1.$
Thus, $l_1$ is a non-zero element of
$I\cap W(n+m)$ and it follows $I=W(n,m)$ as
the above case.

\noindent
If $a_{j_1}=0$ and $j_1=1,$ then from (\ref{s30})
\begin{eqnarray}\label{s45}
& &l_1=c_2 a_{j_2} x^{\beta_1}\partial_1 +\cdots 
+c_t a_{j_t}x^{\beta_t}\partial_1 \\ \nonumber
&+&\cdots \\ \nonumber
&+&c_1 a_{1}x^{\beta_t}\partial_{j_1}
+c_1 b_{11} x^{\beta_1-\epsilon_1}\partial_{j_1}
\\ \nonumber
&+&\cdots \\ \nonumber
&+&c_t a_{1} x^{\beta_t}\partial_{j_t}
+c_t b_{1t} x^{\beta_t-\epsilon_1}\partial_{j_t}.
\end{eqnarray}
Since $l_1\neq 0$, $x^{\beta_1}>x^{\beta_2}$ by assumption.

\noindent
%%%%%%%%%%%%%%%%%%%%%%%%%%Thus, we proved the theorem by lemma 3 and the simplicity
%%%%%%%%%%%%%%%%%%%%%%%%%%of $W(n+m).$

\noindent
We can assume $a_{j_1} = 0.$ 
%%%%%%%%%%%%%%%%%%%and $j_1\neq 1,$ 
Then from (\ref{s30})
%\begin{eqnarray}\label{s50}
%& &l_1=c_2 a_{j_2} x^{\beta_2}\partial_1 +\cdots
%+c_t a_{j_t}x^{\beta_t}\partial_1 \\ \nonumber
%%&+&\cdots \\ \nonumber
%&+&c_1 a_{1}x^{\beta_1}\partial_{j_1} +\cdots
%+c_1 b_{11}x^{\beta_1-\epsilon_1}\partial_{j_1} 
%\\ \nonumber
%&+&\cdots \\ \nonumber
%&+&c_t a_{1}x^{\beta_t}\partial_{j_t}
%+c_t b_{1t}x^{\beta_t-\epsilon_1}\partial_{j_t}.
%\end{eqnarray}
the term $c_1 a_{1}x^{\beta_1}\partial_{j_1}$
is maximal in $l_1,$
so that element of $I\cap W(n+m)$
has a non-zero element. In this case as before $I=W(n,m).$
%%%%%%%%%%%%%%%%$l_1\neq 0$, 
%%%%%%%%%%%%%%%%because of the term $c_1 a_{1}x^{\beta_1}\partial_{j_1}.$ 
%%%%%%%%%%%%%%%%%$c_1 a_1\neq 0$ and
%%%%%%%%%%%%%%%%%$\partial_{j_1}$ is unique in ( ).
%%%%%%%%%%%%%%%%Thus, we have proved the theorem by lemma 3 and the simplicity
%%%%%%%%%%%%%%%%of $W(n+m).$

\noindent
%{\bf Step 2}
%%%%%%%%%%%%%%%%%%%%%%%%Assume $st(l)=2$ for any $l\in I,$  then the theorem is true.
Let $n\geq 2$ and suppose that $I=W(n,m)$ for
an ideal $I$ containing a non-zero element $l$ with
$st(l) \neq n-1.$

\noindent
%{\it Proof.}
%%%%%%%%%%%%%%%%%%%%If $l\in I$ $st(l)=2,$ then $l$ can be written as follows:
Let $l\in I$ and $st(l)=n.$
Then $l$ can be written as follows:
\begin{eqnarray}\label{s55}
& &l=\sum_{\beta,k}c_{\alpha,\beta,k} e^{\alpha }x^{\beta}\partial_k
+\cdots \\ \nonumber
&+&\sum_{\psi,p}c_{\phi,\psi,p} e^{\phi }x^{\psi}\partial_p
\end{eqnarray}
where $\alpha >\cdots >\phi$ and
$\beta >(0,\cdots ,0)$ and $\psi >(0,\cdots ,0)$ by Lemma 1.

%%%%%%%%%%%%%%%%%& &l=\sum_{\beta_i}e^{\alpha x}x^{\beta_i}\partial_k
%%%%%%%%%%%%%%%%%+\sum_{\psi_j}e^{\phi x}x^{\psi_j}\partial_p
If some terms with $e^{0}$ occur in (\ref{s55}), 
then we can assume there is
%%%%%%%%%%%%%%%%%%%%%%%%%$a_u\neq 0$ such that $(a_1,\cdots ,a_u,\cdots ,a_n).$
$a_u\neq 0$ for some $u$ in  $\alpha = (a_1,\cdots ,a_u,\cdots ,a_n).$
%%%%%%%%%%%%%%%%%%Then $$l_2=[\partial_u,[\partial_u,[\cdots ,[\partial_u,l] \cdots]$$
%%%%%%%%%%%%%%%%%%where we applied Lie bracket $(lp(l)+1)$-times; then
Let $$l_2=[\partial_u,[\partial_u,[\cdots ,[\partial_u,l] \cdots]$$
where we applied Lie bracket $(lp(l)+1)$-times. Then
the $(0,\cdots ,0)-$homogeneous part is zero 
, $st(l_2)= n-1$  and  $l_2\neq 0$ by Lemma 3.
Thus, we have $I=W(n,m)$ by induction hypothesis.
%%%%%%%%%%%%%%%%%%%%%%%%of $W(n+m).$

\noindent
%%%%%%%%%%%%%%%%%%%If $\phi\neq (0,\cdots ,0)$, then we have
%%%%%%%%%%%%%%$\phi=(0,\cdots ,0,\phi_j,\cdots ,\phi_u)$ with
%%%%%%%%%%%%%%$\phi_j \neq 0.$ By step 1, we can find
%%%%%%%%%%%%%%$u\in \{1,\cdots ,n+m\}$ such that
If $\phi =(h_1,\cdots ,h_n)\neq (0,\cdots ,0)$, then
we have $h_u \neq 0$ for some $l\leq u \leq n.$
Then we have
\begin{eqnarray}\label{s60}
& &[e^{-\phi}\partial_u,l]=[e^{-\phi }\partial_u \quad , \quad \sum_{\beta,k} c_{\alpha ,\beta,k}e^{\alpha }x^{\beta}\partial_k 
+\cdots + \sum_{\psi,p} c_{\phi ,\psi,p}e^{\phi }x^{\psi}\partial_p ]\\ \nonumber
&=& \sum_{\beta,k} a_u c_{\alpha ,\beta,k}e^{\alpha -\phi }x^{\beta}\partial_k \\ \nonumber
&+&\sum_{\beta,k} b_u c_{\alpha ,\beta,k}e^{\alpha -\phi }x^{\beta-\epsilon_u}\partial_k \\ \nonumber
&+& \sum_{\beta,k} h_k c_{\alpha ,\beta,k}e^{\alpha -\phi }x^{\beta}\partial_u \\ \nonumber
& &\vdots \\ \nonumber
%& &\cdot \\ \nonumber
%& &\cdot \\ \nonumber
&+& \sum_{\psi,p} h_u c_{\phi ,\psi,p}x^{\psi}\partial_p \\ \nonumber
&+& \sum_{\psi,p}\psi_u c_{\phi ,\psi,p}x^{\psi-\epsilon_u}\partial_p \\ \nonumber
&+& \sum_{\psi,p} h_p c_{\phi ,\psi,p}x^{\psi}\partial_u 
\end{eqnarray}
such that $\beta >(0,\cdots ,0), \cdots , \psi >(0,\cdots ,0)$ by 
Lemma 1.

%\begin{eqnarray}\label{s60}
%& &l_3=[e^{-\phi }\partial_k \quad , \quad \sum_{\beta_i} c_{\alpha ,\beta_i}e^{\alpha }x^{\beta_i}\partial_u 
%+ \sum_{\psi_p} c_{\phi ,\psi_p}e^{\phi x}x^{\psi_p}\partial_t ]\\ \nonumber
%&=& \sum a_k c_{\alpha ,\beta_i}e^{\alpha x-\phi x}x^{\beta_i}\partial_u \\ \nonumber
%&+&\sum b_k c_{\alpha ,\beta_i}e^{\alpha x-\phi x}x^{\beta_i}\partial_u \\ \nonumber
%&+& \sum h_u c_{\alpha ,\beta_i}e^{\alpha x-\phi x}x^{\beta_i-\epsilon_u}\partial_k \\ \nonumber
%&+& \sum h_k c_{\phi ,\psi_p}x^{\psi_p}\partial_t \\ \nonumber
%&+& \sum s_k c_{\phi ,\psi_p}x^{\psi_p-\epsilon_k}\partial_t \\ \nonumber
%&+& \sum h_t c_{\phi ,\psi_p}x^{\psi_p}\partial_k .
%\end{eqnarray}

\noindent
$l_3$ is non-zero by the argument similar to the case
$st(l)=1.$
If the $(\alpha -\psi)$-homogeneous component is zero,
then we can apply the induction hypothesis to have $I=W(n,m).$
If the $(\alpha -\psi)$-homogeneous component is non-zero,
then we have
$$l_4=[\partial_j,\cdots ,[\partial_j,l_3]\cdots ]$$
the Lie bracket $(lp(l)+1)$-times.
Then $l_4$ has no $(0,\cdots ,0)-$homogeneous components,
and we can apply induction hypothesis to have $I=W(n,m).$
This completes the proof.
\quad $\Box$

\section{Some Properties of $W(n,m)$.}

Consider the gradations of the Witt algebra $W(n).$ 
The $Z$-gradation of Witt algebra $W(n)$ is well known \cite{Rud}.
In fact $W(n)$ has a  
$Z^k 
-$gradation for $1\leq k \leq n$ as follows:
\begin{eqnarray}\label{s65}
& &W(n)=\bigoplus _{(a_1,\cdots ,a_k)\in Z^k} 
W_{(a_1,\cdots ,a_k)} 
\end{eqnarray}
where
$W_{(a_1,\cdots ,a_k)} $ 
%has a basis
is the subspace of $W(n)$ with basis
$$\{ x_1^{a_1} \cdots x_{u-1}^{a_{u-1}} x_u^{a_u+1 }x_{u+1}^{a_{u+1}} \cdots x_k^{a_k} x_{k+1}^{a_{k+1}} \cdots x_n^{a_n}\partial_u 
\quad |\quad 1\leq u \leq k, 
a_{k+1},\cdots ,a_n\in Z \}\cup$$
$$\{ x_1^{a_1} \cdots x_k^{a_k}x_{k+1}^{a_{k+1}} \cdots x_n^{a_n}\partial_j 
\quad |\quad k+1\leq j \leq n, 
a_{k+1},\cdots ,a_n\in Z \}.$$
%%%%%%%%%%%%%%%%%%%%
%%%%%%%%%%%%%%%%%%%%
%%%%%%%%%%%%%%%%%%%%
%%%%%%%%%%%%%%%%%%%%

Let us denote $W_{(0,\cdots ,0)\in Z^k}$ by $W_0^k,$
where $k\in \{1,\cdots ,n\}.$
Note that $W_0\supset W^2_0 \supset \cdots \supset W^n_0.$

%%%%%%%%%%%%%%%%%%%%\begin{prop}
%%%%%%%%%%%%%%%%%%%%The Lie algebra $W(n)$ is not isomorphic to $W(m),$
%%%%%%%%%%%%%for $n\neq m.$
%%%%%%%%%%%%%\end{prop}
%%%%%%%%%%%%%{\it Proof.}
%%%%%%%%%%%%%If $W(m)$ is isomorphic to $W(n)$, then
%%%%%%%%%%%%%there is an isomorphism $\theta $ from $W(n)$ to $W(m)$.
%%%%%%%%%%%%%Since $W(m)$ has $Z^m-$gradation and
%%%%%%%%%%%%%$W(n)$ has $Z^n-$gradation, thus $\theta (W^n_0)$ 
%%%%%%%%%%%%%%of $W(m)$ 
%%%%%%%%%%%%%%%%%%%%%%%%%%is a subalgebra of $W(m).$ 
%%%%%%%%%%%%%%and $\theta (W^n_0)$ is a 
%%%%%%%%%%%%%%subalgebra of $W(n).$ 
%%%%%%%%%%%%%$W^m_0$ has a basis 
%%%%%%%%%%%%%$\{x_1\partial_1, \cdots , x_m\partial_m \}$ and
%%%%%%%%%%%%%$W^n_0$ has dimension $n.$  
%%%%%%%%%%%%%But the subalgebra of $W(n)$ with basis 
%%%%%%%%%%%%%$\{\theta (x_1\partial_1), \cdots , \theta (x_m\partial_m)\}$ has dimension $m$,
%%%%%%%%%%%%%which contradicts $m\neq n$, therefore the proposition is proved.
%$\theta (W^m_0)$ is not isomorphic to 
%%%%%%%%%%%%%%$W^n_0.$  
%%%%%%%%%%%%%\quad $\Box$
%%%%%%%%%%%%%
%%%%%%%%%%%%%\bigskip

%%%%%%%%%%%%%%%%%%%%%%%%%%%%%%%%%%
%%%%%%%%%%%%%In this section we use the $Z^n-$gradation of $W(n,m)$ for
%%%%%%%%%%%%%the following theorem.
Actually $W(n,m)$ has 
a $Z^{n+m}$-gradation as follows:
$$W(n,m)=\bigoplus
_{(a_1,\cdots, a_n,a_{n+1},\cdots ,a_{n+m})\in Z^{n+m}}
W^{n+m}_{(a_1,\cdots, a_n,a_{n+1},\cdots ,a_{n+m})}$$
where
$W^{n+m}_{(a_1,\cdots, a_n,a_{n+1},\cdots ,a_{n+m})}$ is
a subspace of $W(n,m)$ with basis
\begin{eqnarray*}
& &\{ e^{a_1 x_1} \cdots e^{a_n x_n}x_1^{i_1} \cdots x_n^{i_n} x_{n+1}^{a_{n+1}}\cdots
x_{n+m}^{a_{n+m}}\partial_j\quad |\\
& & 1\leq j\leq n, i_1,\cdots ,i_n\in Z 
\}\\
& &\cup \{e^{a_1 x_1} \cdots e^{a_n x_n}x_1^{p_1} \cdots x_n^{p_n} x_{n+1}^{a_{n+1}}\cdots
x_{n+j}^{a_{n+j}+1} x_{n+j+1}^{a_{n+j+1}} \cdots x_{n+m}^{a_{n+m}}\partial_{n+j}|\\
& & 1\leq j\leq m, p_1,\cdots ,p_n\in Z 
\}.
\end{eqnarray*}

%%%%%%%%%%%%%%%%%%%%%%%%%
%%%%%%%%%%%%%\begin{thm}
%%%%%%%%%%%%%The Lie algebra $W(n,m)$ can not be isomorphic to $W(n_1,m_1),$
%%%%%%%%%%%%%for fixed $n,n_1\in N$ such that $n\neq n_1$ and for
%%%%%%%%%%%%%any $m,m_1\in N.$
%%%%%%%%%%%%%\end{thm}
%%%%%%%%%%%%%{\it Proof.}
%%%%%%%%%%%%%Without loss of generality we can assume $n>n_1.$
%%%%%%%%%%%%%Let $\theta$ be a isomorphism from
%%%%%%%%%%%%%$W(n,m)$ to $W(n_1,m_1).$
%%%%%%%%%%%%%We have the subalgebra $W_0^{n+m}$ of $W(n,m)$ with basis
%%%%%%%%%%%%%$$\{x_1^{i_1} \cdots  x_n^{i_n}\partial_j \quad |\quad j\in \{
%%%%%%%%%%%%%1,\cdots ,n+m\},i_1,\cdots ,i_n\in Z\}.$$ 
%%%%%%%%%%%%%Also we have 
%%%%%%%%%%%%%the subalgebra $W_0^{n_1+m_1}$ of $W(n,m)$ with basis
%%%%%%%%%%%%%$$\{x_1^{j_1} \cdots  x_{n_1}^{j_{n_1}}\partial_k\quad |\quad k\in \{
%%%%%%%%%%%%%1,\cdots ,n_1+m_1\},j_1,\cdots ,j_{n_1}\in Z\}.$$ 
%%%%%%%%%%%%%$W_0^{n+m}/M_1\cong W(n)$  
%%%%%%%%%%%%%where $M_1$ is a maximal abelian
%%%%%%%%%%%%%subalgebra of $W_0^{n+m}$ and
%%%%%%%%%%%%%$W_0^{n_1+m_1 }/M_2 \cong W(n)$
%%%%%%%%%%%%%where $M_2$ is a maximal abelian
%%%%%%%%%%%%%subalgebra of $W_0^{n_1+m_1}.$
%%%%%%%%%%%%%%By above proposition,
%%%%%%%%%%%%%%$W_0^{n+m}$ can not be isomorphic to
%%%%%%%%%%%%%%$W_0^{n_1+m_1 }\supset W(n)$.
%%%%%%%%%%%%%Then $\theta (W(n))\subset W(n_1).$ This contradicts
%%%%%%%%%%%%%the fact that $n_1<n,$ by Proposition 1.
%%%%%%%%%%%%%\noindent
%%%%%%%%%%%%%Thus, there does not exist an isomorphism
%%%%%%%%%%%%%between them.
\quad $\Box$

\bigskip

\noindent
%Let us denote $Z^n$ as the $n$-times of $Z$, i.e.
%$Z \hbox { x } Z \cdots Z \hbox { x } Z.$
%%%%%%%%%%%%%%%%%%%%%%%%%%%%%%%%%%
%Actually $W(n,m)$ has 
%a $Z^{n+m}$-gradation as follows:
%$$W(n,m)=\bigoplus
%_{(a_1,\cdots, a_n,a_{n+1},\cdots ,a_{n+m})\in Z^{n+m}}
%W_{(a_1,\cdots, a_n,a_{n+1},\cdots ,a_{n+m})}$$
%where
%$W_{(a_1,\cdots, a_n,a_{n+1},\cdots ,a_{n+m})}$ is
%%a subspace of $W(n,m)$ with basis
%\begin{eqnarray*}
%& &\{ e^{a_1 x_1} \cdots e^{a_n x_n}x_1^{i_1} \cdots x_n^{i_n} x_{n+1}^{a_{n+1}}\cdots
%x_{n+m}^{a_{n+m}}\partial_j\quad |\\
%& & 1\leq j\leq n, i_1,\cdots ,i_n\in Z 
%\}\\
%%& &\cup \{e^{a_1 x_1} \cdots e^{a_n x_n}x_1^{p_1} \cdots x_n^{p_n} x_{n+1}^{a_{n+1}}\cdots
%x_{n+j}^{a_{n+j}+1} x_{n+j+1}^{a_{n+j+1}} \cdots x_{n+m}^{a_{n+m}}\partial_{n+j}|\\
%& & 1\leq j\leq m, p_1,\cdots ,p_n\in Z 
%\}.
%%\end{eqnarray*}

%%%%%%%%%%%%%%%%%%%%%%%%%%%%%%%%%%%%%%%%%%%%%%%%%%%%%%%%%%%%%%%%%%%%%%%%%%%%%%%%%%%%%%%%%%%%%%%%%%%%%%%
\noindent
$W_{(0,\cdots ,0)}^{n+m}$ is the subalgebra of $W(n,m)$ with
basis
\begin{eqnarray} \label{s70}
& &\{x_1^{i_1}\cdots  x_{n}^{i_{n}} \partial_{j}\quad |\quad 1\leq j \leq n,i_1,\cdots ,i_n\in Z\}
\end{eqnarray}
\begin{eqnarray}\label{s75}
& &\cup \{x_1^{p_1}\cdots  x_{n}^{p_{n}} \partial_{n+j}\quad |\quad 1\leq j \leq m ,p_1,\cdots ,p_n\in Z\}.
\end{eqnarray}
$W_{(0,\cdots ,0)}^{n+m}$ 
can be decomposed as follows:
$$W_{(0,\cdots ,0)}^{n+m}=W(n)\oplus A,$$ 
where $W(n)$ has a basis (\ref{s70}), and $A$ has a basis (\ref{s75}).
$A$ is an abelian ideal of
$W_{(0,\cdots ,0)}^{n+m}.$ 

%\begin{lemma}
%$[W_{(0,\cdots ,0)}^{n+m}, 
%W_{(0,\cdots ,0)}^{n+m}]=A$  and $W^{n+m}_{(0,\cdots ,0)}$ is
%solvable \cite{Hum}.
%\end{lemma}
%{\it Proof.}
%It is straightforward.
%\quad $\Box$
%%\end{document}
%%%%%%%%%%%%%%%%%%%%%%%%%
%\begin{prop}
%$W_{(0,\cdots ,0)}^{n+m}$ is not nilpotent.
%%but it is solvable [1,6]. 
%\end{prop}
%{\it Proof.}
%$$[W_{(0,\cdots ,0)}^{n+m}, 
%[W_{(0,\cdots ,0)}^{n+m},[\cdots  
%[W_{(0,\cdots ,0)}^{n+m},  
%W_{(0,\cdots ,0)}^{n+m}]\cdots  ]=A.$$
%Therefore, 
%$W_{(0,\cdots ,0)}^{n+m}$ is not nilpotent.
%%\noindent
%%From lemma 5, 
%%$W_{(0,\cdots ,0)}^{n+m}$ is 
%%solvable.
%\quad $\Box$
%
%%%%%%%%%%%%%%%%%%%%%%
\bigskip

\noindent
To compare $W(n,m)$ with other algebras, 
we need to find the ad-semisimple
elements. Since $W(n,m)$ has a sub-algebra $W(0,m)$,
we see the element $x^{\epsilon_t}\partial_t,$ where
$t\in \{n+1,\cdots,n+m\}$ and $\epsilon_t=(0,\cdots,0,1,0,\cdots,0),$
is an ad-semisimple element of $W(0,m)$.

\bigskip

\noindent
$W(n,m)$ has a subalgebra $W(n)$ where $n$ comes
from the exponential parts. Thus
$x^{\epsilon_t}\partial_t$ is the candidate of ad-semisimple
elements for $t\in \{1,\cdots,n\}.$ But
$$[e^{-\epsilon_t}\partial_t \quad ,\quad x^{\epsilon_t}\partial_t]
=e^{-\epsilon_t}\partial_t+e^{-\epsilon_t} x^{\epsilon_t}\partial_t.$$ 
This shows $x^{\epsilon_t}\partial_t$ can not 
be an ad semisimple element of $W(n,m)$. 
Thus we have the following proposition.
%%%%%%%%%%%%%%%%%%%%%%%%%%%%%%%%%%%%%%%%%%%%%%%%%%%%%
%%%%%%%%%%%%%%%%%%%%%%%%%%%%%%%%%%%%%%%%%%%%%%%%%%%%%
%%%%%%%%%%%%%%%%%%%%%%%%%%%%%%%%%%%%%%%%%%%%%%%%%%%%%
%%%%%%%%%%%%%%%%%%%%%%%%%%%%%%%%%%%%%%%%%%%%%%%%%%%%%
%%%%%%%%%%%%%%%%%%%%%%%%%%%%%%%%%%%%%%%%%%%%%%%%%%%%%
%%%%%%%%%%%%%%%%%%%%%%%%%%%%%%%%%%%%%%%%%%%%%%%%%%%%%

\begin{prop}
The Lie algebra $W(n,m)$ has an $m-$dimensional maximal
torus which has the basis
$\{x^{\epsilon_t}\partial_t \qquad | \qquad n+1\leq t\leq n+m\}$
with respect to the basis (\ref{s3}).
\end{prop}
{\it Proof.}
Let $l$ be an ad-semisimple element with respect
to the basis (\ref{s3}).
Then $l$ can be written as follows:
$$l=\sum_i f_i\partial_i.$$
Then for any $j$
$$\sum_i \partial_j (f_i)\partial_i=[\partial_j,l]\in F\partial_j.$$
So $\partial_j(f_i)=0$ for $j\neq i$ and $\partial_i(f_i)=a_i\in F.$
It follows that $\partial_j(f_i-a_ix_i)=0$ for all $i,j.$
So $f_i-a_ix_i$ is a constant. Say $f_i-a_ix_i=b_i.$
Thus $f_i=a_ix_i+b_i$ with $a_i,b_i\in F.$

Next 
\begin{eqnarray*}
& &[x_j\partial_j,l]=[x_j\partial_j, \sum_i(a_ix_i+b_i)\partial_i]\\
&=&(a_jx_j-(a_jx_j+b_j))\partial_j=-b_j\partial_j
\end{eqnarray*}
is a multiple of $x_j\partial_j.$ Thus
$b_j=0$ and $l=\sum_i a_i x_i\partial_i.$

Finally if $j\leq n,$ then
\begin{eqnarray*}
& &[e^{x_j}\partial_j,l]=[e^{x_j}\partial_j, \sum_ia_ix_i\partial_i]\\
&=&e^{x_j}(a_j-a_jx_j)\partial_j
\end{eqnarray*}
is a scalar multiple of $e^{x_j}\partial_j$ so $a_j=0.$

Thus $l=\sum_{i=n+1}^{n+m} a_i x_i\partial_i.$
\quad $\Box$

\noindent
Obviously we have following corollary.
\begin{coro}
The Lie algebra $W(n,0)$ has
no ad-diagonal element with respect to the basis
$$\{e^{a_1x_1}\cdots e^{a_nx_n}x_1^{b_1}\cdots x_{n}^{b_n}\partial_i|a_1,\cdots
,a_n,b_1,\cdots ,b_n\in Z,1\leq i\leq n\}.$$
\end{coro}

Consider the Lie subalgebra $W^+(1,0)$ of $W(n,0)$
which is spanned by
$$\{e^{ax}x^i \partial \quad |\quad a,i\in N \}$$
%with the obvious Lie bracket (1). 
This is a subalgebra
of $W(1,0)$ and we have the following proposition.
\begin{prop}
The Lie algebra $W^+(1,0)$ is a non-simple Lie
algebra \cite{Kac}.
\end{prop}
{\it Proof.}
Consider the subalgebra $I_m$ which is spanned by 
$$\{e^{ax}x^i\partial \qquad |a\geq m, i\in N\}$$
where $m\geq 0$ and $a$ is any fixed positive
integer. Then $I_m$ is a non-trivial ideal of $W^+(1,0)$.
\quad $\Box$

\bigskip

\noindent
This proposition shows that if we want a simple 
Lie algebra, we need to include the terms of the
form $e^{ax}x^i \partial$ where $a$ can be any negative integer
%if we want a simple Lie algebra
. We have proved that the 
Lie algebra $W(1,0)$ is simple from Theorem 1.
%%%%%%%%%%%%%%%%%%%%%%%%%%%%%%
%%%%%%%%%%%%%%%%%%%%%%%%%%%%%%
%%%%%%%%%%%%%%%%%%%%%%%%%%%%%%

We have the following Lie embedding:
%We have series of generalized Witt algebras as follows:
$$W(1)\subset W(2)\subset \cdots \subset W(m)\subset \cdots ,$$
$$W(1,1)\subset W(1,2)\subset \cdots \subset W(1,m)\subset \cdots ,$$
$$ \cdots $$
$$W(n,1)\subset W(n,2)\subset \cdots \subset W(n,m)\subset \cdots .$$

\bigskip
\noindent
Let us study a simple Lie algebra which has no toral
element.
Let $0\neq q\in C$, the complex numbers, be a fixed 
non-root of unity $(q^n\neq 1 \hbox { for any } n\in N).$
The skew polynomial ring $C_q[x,y],$ where $yx=qxy,$ has been 
called the quantum plane. The quantum plane can be
localized at the Ore set of powers of $x,y$ to give a non-commutative
associative ring
of Laurent polynomials $C_q[x,y,x^{-1},y^{-1}]$ \cite{Kir}.
\begin{prop}\cite{Kir}
The Lie algebra of derivations of $C_q[x,y,x^{-1},y^{-1}]$
is generated by the inner derivations and the derivations
$D_{\alpha,\beta}(x^iy^j)=(\alpha i +\beta j)x^i y^j.$
\end{prop}
{\it Proof.}
See Theorem 1.2 of \cite{Kir} on page 3757.
\quad $\Box$

Let $V_q$ be the Lie algebra of inner derivations on
$C_q[x,y,x^{-1},y^{-1}].$
Clearly $\{ad(x^h y^k)| (h,k)\in Z \hbox { x } Z-\{(0,0)\} \}$ forms a 
basis of $V_q$, and
$$[ad(x^h y^k),ad(x^ry^s)]=ad([x^hy^k,x^ry^s]).$$

\begin{prop}
The Lie algebra $V_q$ is a simple Lie algebra.
\end{prop}
{\it Proof.}
This case is Theorem 1.3 of \cite{Kir}.
\quad $\Box$

\bigskip

Let us consider the Lie algebra $\bar V_q$ on 
$C_q[x, x^{-1},y,y^{-1}]$ with a different
viewpoint from Kirkman, Processi and Small's idea \cite{Kir}.

The Lie algebra $\bar V_q$ has a basis
$\{x^i y^j \quad | \quad i,j\in Z\}$ and
a Lie bracket on the basis given by
$$[x^i y^j, x^l y^m]= 
x^i y^j x^l y^m
-x^l y^m x^i y^j$$
for any pair of basis elements
$x^i y^j, x^l y^m  \in \bar V_q.$

\begin{coro}
$V_q \cong \bar V_q.$
\end{coro}
{\it Proof.}
It is straightforward.
\quad $\Box$

%%%%%%%%%%%%%%%%%%%%%%%%%%%

\noindent
Consider the Lie algebra $\bar V$ which is spanned by
$$\{(a,i)\quad |\quad a,i\in Z\}$$
with Lie bracket
$$[(a,i),(b,j)]=(q^{bi}-q^{aj})(a+b,i+j).$$

\begin{prop}
$V_q \cong \bar V.$
\end{prop}
{\it Proof.}
Define the $F$-linear map $\theta$ from $V_q$ to $\bar V$
by $\theta (ad (x^h y^k))=(h,k).$
Then clearly $\theta$ is a bijective Lie algebra homomorphism
. Therefore, we have proved the
proposition.
\quad $\Box$

\begin{coro}
$\bar V\cong V_q \cong \bar V_q.$
\end{coro}
{\it Proof.}
It is straightforward.
\quad $\Box$

\bigskip

\noindent
Let us study $\bar V$ more closely.
$(0,0)$ is a basis of the center of $\bar V$
, since $[(0,0),(a,i)]=(0,0)$ for any $(a,i)\in \bar V.$
Let $(h,k)$ be the toral element of $\bar V,$ then for any
$(a,i)\in \bar V,$
$$[(h,k),(a,i)]=(q^{ka}-q^{hi})(a+h,i+k)=f(a,i)$$
for some $f\in C.$
Thus $(h,k)=(0,0).$
%But the element $(0,0)$ is removed for the simplicity of
This contradicts the simplicity of $\bar V\cong V_q.$ 
$(0,0)$
cannot be a basis,
thus, there is no toral element.
The above argument shows that $\bar V$ cannot contain
the Witt algebra $W(1)$ or $W^+(1)$. 
Since $W(n,m)$ contains $W(1)$,    
$W(n,m)$ cannot be isomorphic to $\bar V.$
%$W(n,m)$ contains $W^(1)$.    
%\end{document}

%For $\alpha \in Z^2$:
%define Lie algebras
%$$Y(s): y_{\alpha}=x^{\alpha}D_1 +sx^{\alpha} D_2,$$
%$$V(r,t): v_{\alpha}=(\alpha_2+t)x^{\alpha}D_1 -(\alpha_1+r)x^{\alpha} D_2,$$
%where $r,s,t\in F$ are scalar parameters. 
%The multiplication in each of these algebras is given 
%by (respectively)
%$$[y_{\alpha},y_{\beta}]=((\beta_1-\alpha_1)+s(\beta_2 -\alpha_2))y_{\alpha +\beta},$$
%$$[v_{\alpha},v_{\beta}]=((\alpha_2 \beta_1-\beta_2 \alpha_1)+r(\alpha_2-\beta_2)-t(-\beta_1 +\alpha_1))v_{\alpha +\beta}.$$
%In $V(s)$ let us to find center of this Lie algebra,
%then for all $\beta \in Z^2$ we have,
%$$[y_{\alpha},y_{\beta}]=((\beta_1-\alpha_1)+s(\beta_2 -\alpha_2))y_{\alpha +\beta}=0$$
%iff $\beta_1=\alpha_1$ and $s=0$ or $\beta_2=\alpha_2.$
%Therefore, there is no non-trivial center in $V(s).$

%In $V(r,t)$, we need to find center of $V(r,t),$ for any
%$\beta \in Z^2,$
%$$[v_{\alpha},v_{\beta}]=((\alpha_2 \beta_1-\beta_2 \alpha_1)-t(-\beta_1 +\alpha_1))v_{\alpha +\beta}=0.$$
%This forces
%$$(\alpha_2 \beta_1 -\beta_2 \alpha_1)+r(\alpha_2 -\beta_2)-t(\alpha_1 -\beta_1)=0,$$
%where $r$ and $t$ are fixed non-zero case, there is no center.
%Therefore, we need to check
%$r=t=0$ case.
%\newline
\noindent

\section{Derivations of $W^*(1,0)$}

%Passman's comments and Miller's comments.

In this section we determine all the derivations of
the Lie subalgebra $W^*(1,0)$ of $W^*(n,m).$
Ikeda and Kawamoto found all the derivations of the
Kawamoto algebra $W(G,I)$ in their paper \cite{Ik}.
It is very important to find the derivations of 
a given Lie algebra to compare with other Lie algebras.
Let $L$ be a Lie algebra over any
field $F$. An $F$-linear map $D$ from $L$ to
$L$ is a derivation if 
$D([l_1,l_2])=[D(l_1),l_2]+[l_1,D(l_2)]$
for any $l_1,l_2\in L.$

Let $L$ be a Lie algebra over any
field $F$. 
Define the derivation $D$ of $L$
to be a scalar derivation if 
for all basis elements $l$ of $L,$ 
$D(l)=f_ll$
where $f_l$ is a scalar depending on $l$
%for some scalar $f_l\in F$ 
 \cite{Blo}.

Let $d$ be any additive function of $Z$ to $F.$
Then it is to see that the linear transformation $D_d$
determined by the mapping of the basis elements of 
$W^*(1,0)$ by 
$D_d(e^{ax}x^i\partial)=d(a)e^{ax}x^i\partial$ is a
derivation of 
$W^*(1,0).$  

In this section we consider the stable $F$-algebra
$F[e^{\pm x},x]$ in $F[[x]]$ with $F$-algebra
basis $\{e^{ax}x^i |a\in Z, i\in N\}$ \cite{Kac}.

Let us calculate all the derivations of
the Lie subalgebra
$W^*(1,0)$ of
$W(n,m).$ 
We need the following lemma for the main
theorem.

\begin{lemma}\label{lem100}
Let $D$ be a 
derivation of $W^*(1,0).$ If $D(\partial)=0$, then
$D=f ad_{\partial} + S$ where $f\in F$
and $S$ is a scalar derivation. 
\end{lemma}
{\it Proof.}
Let $D$ be a derivation of
$W^*(1,0)$ such that $D(\partial)=0.$

We have 
\begin{eqnarray}\label{deri1}
& &D([\partial,x\partial])=D(\partial)=0.
\end{eqnarray}

\noindent
On the other hand, this equals
\begin{eqnarray}\label{deri2}
& &D([\partial,x\partial])=[D(\partial),x\partial]+[\partial,D(x\partial)]=[\partial,D(x\partial)].
\end{eqnarray}
%We can put $D(x\partial)=\sum_{a,i} C(a,i,0,1) e^{ax}x^i\partial$ 
We can put $D(x\partial)=f\partial$ 
for $f\in F[e^{\pm x},x].$ 
%using the ordering (\ref{s15}).
%(For the sake of clarity, from now on we generally omit
%subscripts.)
Then from (\ref{deri1}) and (\ref{deri2}),
\begin{eqnarray}\label{deri5}
& &0=[\partial, D(x\partial)] =[\partial,f\partial]=f'\partial.
%\\ \nonumber
%&=&\sum_{a,i} a C(a,i,0,1) e^{ax}x^i\partial
%+\sum_{a,i} i C(a,i,0,1) e^{ax}x^{i-1}\partial
\end{eqnarray}
%for $C(a,i,0,1), C(a,i,0,1)\in F.$
Since $f$ is an analytic function satisfying $f'=0$, we have $f=c$
for $c\in F.$
%This implies $a=0.$
%\noindent
%Set
%$D(x\partial)=\sum C(0,i,0,1) x^i\partial.$
From (\ref{deri5}) we have
$D(x\partial)=c \partial.$

%& &D([x\partial,x^3\partial]) =2D(x^3\partial)=
%\alpha \partial+6C(0,0,0,1)x^2\partial.
%\end{eqnarray}
%On the other hand,
%\begin{eqnarray}\label{deri45}
%& &[D(x\partial),x^3\partial]
%+[x\partial,D(x^3\partial)] \\ \nonumber
%&=&[C(0,0,0,1)\partial
%,x^3\partial]
%+[x\partial
%,D(x^3\partial)] \\ \nonumber
%&=&3C(0,0,0,1)x^2\partial
%+[x\partial
%,\alpha \partial +3C(0,0,0,1)x^2\partial]\\ \nonumber
%&=&3C(0,0,0,1)x^2\partial
%-\alpha \partial +3C(0,0,0,1)x^2\partial.
%\end{eqnarray}
%This implies $\alpha =0.$
%

Inductively, we assume
\begin{eqnarray}\label{deri50}
& &D(x^n\partial)=ncx^{n-1}\partial
\end{eqnarray}
for $n\in N.$

\noindent
We have
\begin{eqnarray}\label{deri60}
& &D([\partial,x^{n+1}\partial]) =(n+1)D(x^n\partial)
=n(n+1)c x^{n-1}\partial.
\end{eqnarray}
On the other hand, this equals
\begin{eqnarray}\label{deri65}
& &D([\partial,x^{n+1}\partial])
=[D(\partial),x^{n+1}\partial]
+[\partial,D(x^{n+1}\partial)]
=[\partial,D(x^{n+1}\partial)]. 
\end{eqnarray}
From (\ref{deri60}) and (\ref{deri65}),
we have
$$D(x^{n+1}\partial)=\phi \partial+(n+1)c x^n\partial$$
for some $\phi \in F.$

\noindent
Also we have
\begin{eqnarray}\label{deri70}
& &D([x\partial,x^{n+1}\partial]) =nD(x^{n+1}\partial)
=n\phi \partial
+n(n+1) c x^{n}\partial.
\end{eqnarray}
On the other hand, this equals
\begin{eqnarray}\label{deri75}
& &D([x\partial,x^{n+1}\partial])
=[D(x\partial),x^{n+1}\partial]
+[x\partial,D(x^{n+1}\partial)] \\ \nonumber
&=&c[\partial, x^{n+1}\partial]
+[x\partial, 
\phi \partial +(n+1)c x^n\partial ]\\ \nonumber
&=&c(n+1)x^n\partial-\phi \partial
+(n-1)(n+1)c x^n\partial \\ \nonumber
&=&c(n+1)nx^n\partial-\phi \partial.
\end{eqnarray}
From (\ref{deri70}) and (\ref{deri75}), we have
$\phi =0.$

Therefore, for any $n\in N$
%(\ref{deri50}) holds.
we have $D(x^n\partial)=ncx^{n-1}\partial.$

%\noindent
%Now $D(e^x\partial)$ can be written as a sum of
%different homogeneous components as follows:
%\begin{eqnarray}\label{deri77}
%& & D(e^x\partial)=C(b,j,1,0) e^{bx} x^j\partial +\cdots +
 %C(b,0,1,0) e^{bx} \partial 
%\\ \nonumber
%&+&\cdots \\ \nonumber
%&+&C(d,l,1,0) e^{dx} x^l\partial +\cdots 
%+C(d,0,1,0) e^{dx} x^l\partial,
%\end{eqnarray}
%%where $*$ is a 
%%remaining $b$-homogeneous and $e^{bx}x^j\partial$ is a maximal
%using the ordering (\ref{s9}) and the gradation (\ref{s15}) such that
%%basis element in $D(e^x\partial)$ whose
%%coefficient is non-zero 
%$b >\cdots >d$ and \\ 
%$C(b,j,1,0),\cdots ,C(d,0,1,0)\in F$
%with $C(b,j,1,0)\neq 0, \cdots ,
%C(d,l,1,0)\neq 0.$

We have
\begin{eqnarray}\label{deri80}
& &D([\partial,e^{x}\partial]) =D(e^x\partial)
.
\end{eqnarray}
On the other hand, this equals
\begin{eqnarray}\label{deri85}
& &D([\partial,e^x\partial])=[D(\partial),e^{x}\partial]
+[\partial,D(e^{x}\partial)]
=[\partial,D(e^{x}\partial)].
\end{eqnarray}
From (\ref{deri80})
and (\ref{deri85}), we have 
$$\partial (D(e^x \partial ))=  D(e^x\partial).$$
Thus, we have
$$D(e^x \partial )= d e^x\partial$$
for $d \in F.$

Inductively assume
$D(e^{px}\partial)=pde^{px}\partial $ for $1\neq p\in N.$
%Let us show by induction that
%$$D(e^{px} x^{n}\partial)=dp e^{px} x^{n} \partial
%+nc e^{px} x^{n-1} \partial$$ for $p\in Z$ and $n\in N.$
We have
\begin{eqnarray}\label{deri160}
& &D([e^x\partial, e^{px}\partial])=(p-1)D(e^{px+x}\partial).
\end{eqnarray}
On the other hand, this equals
\begin{eqnarray}\label{deri165}
& &[D(e^x\partial), e^{px}\partial]
+[e^x\partial, D(e^{px}\partial)] \\ \nonumber
&=&[de^x\partial, e^{px}\partial]
+[e^x\partial, pde^{px}\partial]\\ \nonumber
&=&(p-1)(p+1)d e^{px+x}\partial.
\end{eqnarray}
From (\ref{deri160}) and (\ref{deri165}), we have
$$D(e^{px+x} \partial)=d(p+1) e^{px+x} \partial.$$
Therefore, we have
$D(e^{px}\partial)=pde^{px}\partial $ for $p\in N.$
                                                    %%%%%%%%%%%%%%%%%%%%%%%%%%%%%%%%
%%%%%%%%%%%%%%%%%%%%%%%%%%%%%%%%

%%%%%%%%%%%%%%%%%%%%%%%%%%%%%%%%%%%%%%%%%%%%%%%%%%%%%%%%%%%%%%%%%%%%%%%%%%%%%%%%%
%%%%%%%%%%%%%%%%%%%%%%%%%%%%%%%%%%%%%%%%%%%%%%%%%%%%%%%%%%%%%%%%%%%%%%%%%%%%%%%%%
%%%%%%%%%%%%%%%%%%%%%%%%%%%%%%%%%%%%%%%%%%%%%%%%%%%%%%%%%%%%%%%%%%%%%%%%%%%%%%%%%
%%%%%%%%%%%%%%%%%%%%%%%%%%%%%%%%%%%%%%%%%%%%%%%%%%%%%%%%%%%%%%%%%%%%%%%%%%%%%%%%%
%%%%%%%%%%%%%%%%%%%%%%%%%%%%%%%%%%%%%%%%%%%%%%%%%%%%%%%%%%%%%%%%%%%%%%%%%%%%%%%%%
Also we have
\begin{eqnarray}\label{deri82}
& &D([\partial,e^{-x}\partial]) =-D(e^{-x}\partial)
.
\end{eqnarray}
On the other hand, this equals
\begin{eqnarray}\label{deri84}
& &D([\partial,e^{-x}\partial])=[D(\partial),e^{-x}\partial]
+[\partial,D(e^{-x}\partial)]
=[\partial,D(e^{-x}\partial)].
\end{eqnarray}
From (\ref{deri82})
and (\ref{deri84}), we have 
$$\partial (D(e^{-x} \partial ))=  D(e^{-x}\partial).$$
Thus, we have
$$D(e^{-x} \partial )= h e^{-x}\partial$$
for $h \in F.$

%%%%%%%%%%%%%%%%%%%
%%%%%%%%%%%%%%%%%%%
%%%%%%%%%%%%%%%%%%%
%%%%%%%%%%%%%%%%%%%
%%%%%%%%%%%%%%%%%%%
%%%%%%%%%%%%%%%%%%%
%%%%%%%%%%%%%%%%%%%
We have 
\begin{eqnarray}\label{deri100}
& &D([e^{-x}\partial, e^x\partial])=2D(\partial)=0.
\end{eqnarray}
On the other hand, this equals
\begin{eqnarray}\label{deri110}
& &D([e^{-x}\partial,e^{x}\partial])=[D(e^{-x}\partial),e^{x}\partial]
+[e^{-x}\partial,D(e^{x}\partial)]\\ \nonumber
&=&[he^{-x}\partial,e^{x}\partial]+[e^{-x}\partial,de^{x}\partial].
\end{eqnarray}
From (\ref{deri100})
and (\ref{deri110}), we have 
$h=-d.$

Inductively assume
$D(e^{px}\partial)=pde^{px}\partial $ where $p$ is a fixed
negative integer.
We have
\begin{eqnarray}\label{deri167}
& &D([e^{-x}\partial, e^{px}\partial])=(p+1)D(e^{px-x}\partial).
\end{eqnarray}
On the other hand, this equals
\begin{eqnarray}\label{deri168}
& &[D(e^{-x}\partial), e^{px}\partial]
+[e^{-x}\partial, D(e^{px}\partial)] \\ \nonumber
&=&[-de^{-x}\partial, e^{px}\partial]
+[e^{-x}\partial, pde^{px}\partial]\\ \nonumber
&=&(p-1)(p+1)d e^{px-x}\partial.
\end{eqnarray}
From (\ref{deri167}) and (\ref{deri168}), we have
$D(e^{px}\partial)=pde^{px}\partial $ holds for any negative
integer $p.$
Therefore,
$D(e^{px}\partial)=pde^{px}\partial $ holds for any $p\in Z.$

%%%%%%%%%%%%%%%%%%%%%%%%%%%%%%%%
%%%%%%%%%%%%%%%%%%%%%%%%%%%%%%%%
%%%%%%%%%%%%%%%%%%%%%%%%%%%%%%%%
%%%%%%%%%%%%%%%%%%%%%%%%%%%%%%%%
%%%%%%%%%%%%%%%%5

Let us show by induction that
$$D(e^{px} x^{n}\partial)=dp e^{px} x^{n} \partial
+nc e^{px} x^{n-1} \partial$$ for $p\in Z$ and $n\in N.$
%\noindent
%Also we have
%\begin{eqnarray}\label{deri160}
%& &D([e^x\partial,e^{x}x^{n+1}\partial]) =(n+1)D(e^{2x} x^{n}\partial)
%%=(n+1)D(e^{2x}x^n\partial).
%.
%\end{eqnarray}
%On the other hand,
%\begin{eqnarray}\label{deri165}
%& &[D(e^x\partial),e^{x}x^{n+1}\partial]
%+[e^x\partial,D(x^{n+1}e^{x}\partial)] \\ \nonumber
%%&=&[\alpha e^x\partial,e^x x^{n+1}\partial]+
%[e^x\partial,(n+1)C(0,0,0,1)x^{n}e^{x}\partial+\alpha e^x x^{n+1}\partial] \\ \nonumber
%&=&\alpha (n+1)e^{2x} x^{n}\partial
%+n (n+1)C(0,0,0,1)e^{2x} x^{n-1}\partial
%+\alpha (n+1)e^{2x} x^{n}\partial.
%\end{eqnarray}
%This gives
%
%$$D(e^{2x} x^{n}\partial)=\alpha x^{n} e^x\partial
%+nC(0,0,0,1) e^{2x}x^{n-1}\partial.$$ 

Inductively assume
$$D(e^{px} x^{n}\partial)=pd e^{px}x^n\partial
+nc e^{px}x^{n-1}\partial$$ for any $p\in Z$ and a fixed $n\in N.$
We have
\begin{eqnarray}\label{deri170}
& &D([x\partial,e^{px}x^{n}\partial]) =pD(e^{px} x^{n+1}\partial)
+(n-1)D(e^{px}x^n\partial).
\end{eqnarray}
On the other hand, this equals
\begin{eqnarray}\label{deri175}
& &[D(x\partial),e^{px}x^{n}\partial]
+[x\partial,D(e^{px}x^n\partial)].
\end{eqnarray}
From (\ref{deri170})
and (\ref{deri175}), we have
\begin{eqnarray}\label{deri185}
& &pD(e^{px}x^{n+1}\partial)=
-(n-1)pd e^{px}x^{n}\partial
-(n-1)nc e^{px}x^{n-1}\partial \\ \nonumber
&+&pce^{px}x^{n}\partial
+nce^{px}x^{n-1}\partial\\ \nonumber
&+&[x\partial, pd e^{px}x^{n}\partial
+nc e^{px}x^{n-1}\partial ]\\ \nonumber
&=&-(n-1)pd e^{px}x^n\partial
-n(n-1)c e^{px}x^{n-1}\partial
+pc e^{px}x^n\partial \\ \nonumber
&+&nc e^{px}x^{n-1}\partial 
+p^2d e^{px}x^{n+1}\partial
+p(n-1)d e^{px}x^{n}\partial\\ \nonumber
&+&pnc e^{px}x^{n}\partial
+nc(n-2) e^{px}x^{n-1}\partial.
\end{eqnarray}
Therefore, we have
$$D(e^{px} x^{n+1}\partial)=pd e^{px}x^{n+1} \partial
+(n+1)c e^{px}x^{n}\partial$$ for any $p\in Z$ and any $n\in N.$

Therefore, $D=c ad_{\partial} + fS$
for an appropriate $f\in F.$
\quad $\Box$

%\begin{theorem}
\begin{thm}
Each derivation of $W^*(1,0)$ can be written
as a sum of an inner derivation
and a scalar derivation.
\end{thm}
%\end{theorem}
{\it Proof.}
%The results from the fact that
%for any derivation $D$ of $W^*(1,0),$
%$D(\partial)$ can be written as the sum of various
%homogeneous elements, and all non-zero homogeneous
%elements are the sum of inner derivations
%and the zero component is a sum of a scalar derivation and 
%inner
%derivations
Let $D$ be any derivation of $W^*(1,0)$. Then
$D(\partial)=f\partial$ for some $f\in F[e^{\pm x},x].$
Since $\partial :
 F[e^{\pm x},x]\to 
 F[e^{\pm x},x]
$ is onto map, there is a function
$g\in
 F[e^{\pm x},x]$ such that $\partial(g)=f.$
Then $ad_{\partial}(g)=[\partial, g\partial]=f\partial=D(\partial).$
We have $(D-ad_{g\partial})(\partial)=0.$
By Lemma
\ref{lem100} we have 
%$D-ad_{\partial}=c ad_{\partial}+S.$
$D=ad_{g\partial}+c ad_{\partial}+S.$
Therefore, we have proven the theorem.
\quad $\Box$

%\end{document}

\section{Some comments on the generalized Witt algebra}

Let $F[x_1,x_2]$ be a polynomial ring. Consider the derivations
$(x_1 x_2\partial_1)$ and
$(x_1 x_2\partial_2)$ on $F[x_1,x_2]$
such that
$$(x_1 x_2\partial_1)(x_1^i x_2^j)=
x_1^i x_2^{1+j},$$
for any 
$x_1^i x_2^j\in F[x_1,x_2]$ and similarly 
%we have
%the above calculation 
for $(x_1 x_2\partial_2).$ 
%on $F[x_1,x_2].$

\noindent
Consider the vector space $L$ over $F$ with basis 
$$\{x_1^i x_2^j (x_1 x_2\partial_1),
x_1^a x_2^b (x_1 x_2\partial_2)| a,b,i,j\in N\}.$$

Then we can define a bracket on the basis as follows [compare with (1)]:
$$[x_1^{i_1} x_2^{j_1}(x_1 x_2\partial_1),
x_1^{i_2} x_2^{j_2}(x_1 x_2\partial_1)]
=(i_2-i_1)x_1^{i_1+i_2} x_2^{j_1+j_2+1}(x_1 x_2\partial_1)$$
and
\begin{eqnarray*}
& &[x_1^{i_1} x_2^{j_1}(x_1 x_2\partial_1),
x_1^{a_1} x_2^{b_2}(x_1 x_2\partial_2)]\\
&=&a_1 x_1^{i_1+a_1} x_2^{j_1+b_2+1}(x_1 x_2\partial_2)
-j_1 x_1^{i_1+a_1+1} x_2^{j_1+b_2}(x_1 x_2\partial_1)
\end{eqnarray*}
and
$$[x_1^{a_1} x_2^{b_1}(x_1 x_2\partial_2),
x_1^{a_2} x_2^{b_2}(x_1 x_2\partial_2)]
=(b_2-b_1)x_1^{i_1+i_2+1} x_2^{j_1+j_2}(x_1 x_2\partial_2),$$
and we extend linearly to $L,$
\noindent
%Then the Jacobi identity does not hold. 
but this bracket does not satisfy the Jacobi identity,
since $(x_1x_2\partial_1)$ and $(x_1x_2\partial_2$
do not commute as derivations.
Thus, it is important to find the generalized derivations
on the polynomial ring to make a Lie algebra.
Please refer to the paper \cite{Pas} for the Witt-type
algebra using generalized derivations on an associative,
commutative $K$-algebra where $K$ is a field of any 
characteristic.

%Thus, $L$ is not a 
%Lie algebra. We should be careful on the action on 
%a basis for the Lie bracket.

\bigskip

\noindent
Following is a quite interesting conjecture on the Witt algebra
\cite{Kawa},\cite{Rud}.

\noindent
{\bf Conjecture.}
If $\theta$ is any automorphism of the Witt algebra $W^+(n),$
then $\theta (x_i\partial_i)=a x_j\partial_j +b\partial_j,$
where $0\neq a\in F$ and $b\in F.$

\bigskip

\noindent
It is not very difficult to prove this conjecture
on $W^+(2)$  and $W^+(1).$

Following is a quite interesting conjecture on the 
derivations of $W(n,m)$ and its subalgebras.

\bigskip

\noindent
{\bf Conjecture.}
Each derivation of $W(n,m)$ can be written as a sum
of an inner derivation and a scalar derivation.

% This is uwrefer.tex

\end{document}